%%%%%%%%%%%%%%%%%%%%%%%%%%%%%%%%%%%%%%%%%%%%%%%%%%%%
% typoref.tex. V : January 18, 2000.
% Author : Anthony PHAN
% Warning : syntaxe +- LaTeX
% Sources :
% T. Lachand--Robert, ``La Ma\^\i trise de \TeX'',
% R\'ef\'erences crois\'ees;
% latex.ltx's sources;
% and of course the \TeX book.
%%%%%%%%%%%%%%%%%%%%%%%%%%%%%%%%%%%%%%%%%%%%%%%%%%%%%
%
\catcode`@=11
%
% style (look at the behavior of \item dans \bibitem too,
% and at one ,\  in \re@dreferenceslist)
% Feel free to change: 	\bibn@me (title like ``R\'ef\'erences'')
%			\bibliographym@rk (general style)
%
\def\bibn@me{R\'ef\'erences}
\def\bibliographym@rk{\centerline{{\sc\bibn@me}}
	\sectionmark\section{\ignorespaces}{\unskip\bibn@me}
	\bigbreak\bgroup
	\ifx\ninepoint\undefined\relax\else\ninepoint\fi}
%
% Beware of the \bgroup: it will be closed by \endthebibliography
%
% \refsp@ce is the spacing command that appens between multiple
% references.
%
\let\refsp@ce=\
\let\bibleftm@rk=[
\let\bibrightm@rk=]
%
% if you want more space between brackets...
%\let\refsp@ce=\thinspace
%\def\bibleftm@rk{[\thinspace}
%\def\bibrightm@rk{\thinspace]}
%
% frenchy stuff
%
\def\numero{n\raise.82ex\hbox{$\fam0\scriptscriptstyle o$}~\ignorespaces}
%
% new variables
%
\newcount\equationc@unt
\newcount\bibc@unt
\newif\ifref@changes\ref@changesfalse
\newif\ifpageref@changes\ref@changesfalse
\newif\ifbib@changes\bib@changesfalse
\newif\ifref@undefined\ref@undefinedfalse
\newif\ifpageref@undefined\ref@undefinedfalse
\newif\ifbib@undefined\bib@undefinedfalse
\newwrite\@auxout
%
% mark an equation
%
\def\eqnum{\global\advance\equationc@unt by 1%
\edef\lastref{\number\equationc@unt}%
\eqno{(\lastref)}}
%
% One can reference anything, just copy the former macro
% and use it so: \machin \label{truc}
% In machin you would have defined \lastref by some number
% or any text.
%
% References macros
%
% The next macros are the core of \ref and \cite commands.
% Its first argument may be ref, pageref or bib.
%
% It is too tricky to be explained.
% (It is a bit recursive.)
% It allows using \cite or \ref or ...
% with arbitrary many arguments,
% for instance:
% \cite{knuth1,knuth2,ma pomme}
%
% First argument is always ref, pageref or bib.
%
\def\re@dreferences#1#2{{%
	\re@dreferenceslist{#1}#2,\undefined\@@}}
\def\re@dreferenceslist#1#2,#3\@@{\def\next{#2}%
	\expandafter\ifx\csname#1@@\meaning\next\endcsname\relax
	??\immediate\write16
	{Warning, #1-reference "\next" on page \the\pageno\space
	is undefined.}%
	\global\csname#1@undefinedtrue\endcsname
	\else\csname#1@@\meaning\next\endcsname\fi
	\ifx#3\undefined\relax
	\else,\refsp@ce\re@dreferenceslist{#1}#3\@@\fi}
%
% notice that the former ``,\refsp@ce'' will separate
% multiple arguments. But beware of spaces
% while defining a reference or calling for it!
%
% tricky thing: \newlabel has two arguments
% {labelname}{{\lastref}{\pageref}}
% The second argument is read as two arguments
% by \newl@bel. This was necessary to get
% a jobname.aux containing the same syntax
% LaTeX would produce and use.
%
\def\newlabel#1#2{{\def\next{#1}\newl@bel#2}}
\def\newl@bel#1#2{%
	\expandafter\xdef\csname ref@@\meaning\next\endcsname{#1}%
	\expandafter\xdef\csname pageref@@\meaning\next\endcsname{#2}}
\def\label#1{{%
	\toks0={#1}\message{ref(\lastref) \the\toks0,}%
	\ignorespaces\immediate\write\@auxout%
	{\noexpand\newlabel{\the\toks0}{{\lastref}{\the\pageno}}}%
	\def\next{#1}%
	\expandafter\ifx\csname ref@@\meaning\next\endcsname\lastref%
	\else\global\ref@changestrue\fi%
	\newlabel{#1}{{\lastref}{\the\pageno}}}}
\def\ref#1{\re@dreferences{ref}{#1}}
\def\pageref#1{\re@dreferences{pageref}{#1}}
%
% bibliography macros
%
\def\bibcite#1#2{{\def\next{#1}%
	\expandafter\xdef\csname bib@@\meaning\next\endcsname{#2}}}
\def\cite#1{\bibleftm@rk\re@dreferences{bib}{#1}\bibrightm@rk}
%
% The argument of \beginthebibliography
% is any sequence of numerals which will represent
% the maximum \item's length. If you have less than 9
% \bibitem's, this argument may be {any numeral}.
% if you have between 100 and 999 \bibitem's
% this argument may be {any three numerals},
% and so on.
%
\def\beginthebibliography#1{\bibliographym@rk
	\setbox0\hbox{\bibleftm@rk#1\bibrightm@rk\enspace}
	\parindent=\wd0
	\global\bibc@unt=0
	\def\bibitem##1{\global\advance\bibc@unt by 1
		\edef\lastref{\number\bibc@unt}
		{\toks0={##1}
		\message{bib[\lastref] \the\toks0,}%
		\immediate\write\@auxout
		{\noexpand\bibcite{\the\toks0}{\lastref}}}
		\def\next{##1}%
		\expandafter\ifx
		\csname bib@@\meaning\next\endcsname\lastref
		\else\global\bib@changestrue\fi%
		\bibcite{##1}{\lastref}
		\medbreak
		\item{\hfill\bibleftm@rk\lastref\bibrightm@rk}%
		}
	}
\def\endthebibliography{\egroup\par}
%
% THE NEXT MACRO MUST BE INCLUDED
% IN THE \BYE COMMAND. FOR INSTANCE:
%
% \catcode`@=11
% \outer\def\bye{\@closeaux
% 	\par\vfill\supereject\end}
% \catcode`@=12
%
\def\@closeaux{\closeout\@auxout
	\ifref@changes\immediate\write16
	{Warning, changes in references.}\fi
	\ifpageref@changes\immediate\write16
	{Warning, changes in page references.}\fi
	\ifbib@changes\immediate\write16
	{Warning, changes in bibliography.}\fi
	\ifref@undefined\immediate\write16
	{Warning, references undefined.}\fi
	\ifpageref@undefined\immediate\write16
	{Warning, page references undefined.}\fi
	\ifbib@undefined\immediate\write16
	{Warning, citations undefined.}\fi}
%
% initialization of jobname.aux
%
\immediate\openin\@auxout=\jobname.aux
\ifeof\@auxout \immediate\write16
  {Creating file \jobname.aux}
\immediate\closein\@auxout
\immediate\openout\@auxout=\jobname.aux
\immediate\write\@auxout {\relax}%
\immediate\closeout\@auxout
\else\immediate\closein\@auxout\fi
%
% Let's read this file and open it out
%
\input\jobname.aux
\immediate\openout\@auxout=\jobname.aux
% this file will be closed by \bye.
%
% That's all, folks!
%

\def\bibn@me{R\'ef\'erences bibliographiques}
%\input typpo
%
%\catcode`@=11
\def\bibliographym@rk{\bgroup}
%
% \bye est modifie pour la biblio et la table des matieres
%
\outer\def\bye{ 	\par\vfill\supereject\end}

\def\R{{\bf R}}

\overfullrule=0pt

\magnification=1200

  \def\pro{\noindent {\bf{Proof : }}}

\def\house#1{\setbox1=\hbox{$\,#1\,$}%
\dimen1=\ht1 \advance\dimen1 by 2pt \dimen2=\dp1 \advance\dimen2 by 2pt
\setbox1=\hbox{\vrule height\dimen1 depth\dimen2\box1\vrule}%
\setbox1=\vbox{\hrule\box1}%
\advance\dimen1 by .4pt \ht1=\dimen1
\advance\dimen2 by .4pt \dp1=\dimen2 \box1\relax}

 \def\la{{\lambda}} \def\deg{{\rm deg}}
   
  \def\eps{{\varepsilon}}

  \def\noi{\noindent}

\def\build#1_#2^#3{\mathrel{\mathop{\kern 0pt#1}\limits_{#2}^{#3}}}

\def\date {le\ {\the\day}\ \ifcase\month\or
janvier\or fevrier\or mars\or avril\or mai\or juin\or juillet\or
ao\^ut\or septembre\or octobre\or novembre\or
d\'ecembre\fi\ {\oldstyle\the\year}}

\font\fivegoth=eufm5 \font\sevengoth=eufm7 \font\tengoth=eufm10

\newfam\gothfam \scriptscriptfont\gothfam=\fivegoth
\textfont\gothfam=\tengoth \scriptfont\gothfam=\sevengoth

\def\cqfd{\unskip\kern 6pt\penalty 500 \raise 0pt\hbox{\vrule\vbox
to6pt{\hrule width 6pt \vfill\hrule}\vrule}\par}

\def\pro{\noindent {\it Proof. }}

\def\smallsquare{\vbox{\hrule\hbox{\vrule height 1 ex\kern 1 ex\vrule}\hrule}}
\def\cqfd{\hfill \smallsquare\vskip 3mm}

\def\hla{{\widehat \lambda}}

\def\RR{{\bf R}}

\def\vp{{\bf p}}
\def\vy{{\bf y}}
\def\ZZ{{\bf Z}}

\def\QQ{{\bf Q}}

%%%%%%%%%%%%%%%%%%%%%%%%%%%%%%%%%%%%%%%%%%%%%%%%%%%%%%%%%%%%%%%%%%%
\vskip 2mm

\centerline{\bf On simultaneous rational approximation
to a real number}
\smallskip

\centerline{\bf and its integral powers, II}

\vskip 13mm

\centerline{Dmitry B{\sevenrm ADZIAHIN} and Yann B{\sevenrm UGEAUD}  \footnote{}{\rm
}}

{\narrower\narrower
\vskip 15mm

\proclaim Abstract.
{For a positive integer $n$ and a
real number $\xi$, let $\lambda_n (\xi)$ denote the supremum
of the real numbers $\lambda$ for which there are arbitrarily
large positive integers $q$ such that
$|| q \xi ||, || q \xi^2 ||, \ldots , ||q \xi^n||$ are all
less than $q^{-\lambda}$. Here, $|| \cdot ||$ denotes the
distance to the nearest integer. We establish new results on the Hausdorff dimension
of the set of real numbers $\xi$ such that $\lambda_n (\xi)$
is equal (or greater than or equal) to a given value.}

}

\vskip 12mm

\centerline{\bf 1. Introduction}

\vskip 6mm

In 1932, in order to define
his classification of real numbers, Mahler \cite{Mah32}
introduced the exponents of Diophantine
approximation $w_n$, which measure how small an integer linear form in the first $n$
powers of a given real number can be.

\proclaim Definition 1.1.
Let $n \ge 1$ be an integer and $\xi$ a real number.
We denote by $w_n(\xi)$
the supremum of the real numbers $w$ such that, for
arbitrarily large real numbers $X$, the inequalities
$$
0 < |x_n \xi^n + \ldots + x_1 \xi + x_0|
\le X^{-w}, \qquad  \max_{0
\le m \le n} \, |x_m| \le X,
$$
have a solution in integers $x_0, \ldots, x_n$.

%The Dirichlet theorem implies that $w_n(\xi)$
%is at least equal to $n$ for every
%real number $\xi$ which is not algebraic of degree at most $n$.
%Sprind\v zuk \cite{Spr69} showed that there is equality for almost all $\xi$, with
%respect to the Lebesgue measure. Furthermore, it follows from the Schmidt Subspace Theorem
%that $w_n(\xi) = \min\{n, d-1\}$ for every positive
%integer $n$ and every real algebraic number $\xi$
%of degree $d$.
We refer to \cite{BuLiv,BuDurham}
for an overview of the known results on the
exponents $w_n$. In particular, it follows from the Schmidt Subspace Theorem
that $w_n(\xi) = \min\{n, d-1\}$ for every positive
integer $n$ and every real algebraic number $\xi$ of degree $d$.
In the sequel, by spectrum of a function,
we mean the set of values
taken by this function at transcendental real numbers.

It is easy to apply the theory of continued fractions
to show that the spectrum of $w_1$ is equal to the
whole interval $[1, + \infty]$. Moreover, the classical
Jarn\'\i k--Besicovich theorem \cite{Jar31} asserts that, for any $w \ge 1$, we have
$$
\dim \{ \xi \in \R : w_1 (\xi) \ge w \} =
\dim \{ \xi \in \R : w_1 (\xi) = w \} =
{2 \over 1 + w}.    \eqno (1.1)
$$
Here, and throughout this paper, $1 / + \infty$ is understood to be $0$ and $\dim$ stands for the
Hausdorff dimension.
To be precise, the Jarn\'\i k--Besicovich theorem concerns the set   %%y
$\{ \xi \in \R : w_1 (\xi) \ge w \}$ and not the level set $\{ \xi \in \R : w_1 (\xi) = w \}$.
However, we easily deduce (1.1) from \cite{Jar31}. In the sequel, we state several metric results
on level sets which, sometimes, are not explicitly stated in the original papers,
but whose validity is known. %%y
For $n \ge 2$,
the fact that the spectrum of $w_n$ equals $[n, + \infty]$
is an immediate consequence of the extension of
(1.1) established in 1983 by Bernik \cite{Ber},  which states that
$$
\dim \{ \xi \in \R : w_n (\xi) \ge w \} = {n + 1 \over w + 1},  \eqno (1.2)
$$
for every positive integer $n$ and every real number $w$ with $w \ge n$.

Another exponent of Diophantine approximation,
introduced in \cite{BuLa05a}, measures the
quality of the simultaneous rational approximation to
the first $n$ integral powers of a real number by rational
numbers with the same denominator.

\proclaim Definition 1.2.
Let $n \ge 1$ be an integer and $\xi$ a real number.
We denote by $\la_n(\xi)$
the supremum  of the real numbers
$\lambda$ such that, for
arbitrarily large real numbers $X$, the inequalities
$$
0 < |x_0| \le X, \qquad \max_{1 \le m \le n} \,
|x_0 \xi^m - x_m| \le X^{-\lambda},     \eqno (1.3)
$$
have a solution in integers $x_0, \ldots, x_n$.

Observe that $\lambda_1$ and $w_1$ coincide.
The Dirichlet theorem implies that $\lambda_n(\xi)$
is at least equal to $1/n$ for every
real number $\xi$ which is not algebraic
of degree at most $n$.
%The combination of
%Sprind\v zuk's above mentioned result with a
%classical transference principle shows that
Furthermore, there is equality for almost all $\xi$, with
respect to the Lebesgue measure; see \cite{Bu10,BuDurham,Schl16} for
further results.
%Furthermore, it follows from the Schmidt Subspace Theorem
%that $\lambda_n(\xi) = \max\{1/n, 1/(d-1)\}$ for every positive
%integer $n$ and every real algebraic number $\xi$
%of degree $d$.
The following question %extends Problem 5.5 from \cite{BuLa07}.
reproduces Problems 2.9 and 2.10 of the survey \cite{BuFARF}.

\proclaim Problem 1.3.
Let $n \ge 1$ be an integer. Is the
spectrum of the function $\lambda_n$ equal to $[1/n,+\infty]$?
For $\lambda \ge 1/n$, what are the Hausdorff dimensions
of the set $\{\xi \in \R : \lambda_n (\xi) \ge \lambda\}$
and of the level set $\{\xi \in \R : \lambda_n (\xi) = \lambda\}$ ?

The above mentioned Jarn\'\i k--Besicovich theorem
answers the case $n=1$ of Problem 1.3. For $n \ge 2$, the state-of-the-art is as follows.
It has been proved in \cite{Bu10} that, for any positive integer $n$ and any
real number $\lambda$ with $\lambda \ge 1$, we can construct
explicitly uncountably many real numbers $\xi$
such that $\lambda_n (\xi) = \lambda$.
Since any real number $\xi$ such that $w_1 (\xi)$ is infinite satisfies
$\lambda_n (\xi) = + \infty$ (see Corollary 3.2 of \cite{Bu10}), we get that
the spectrum of $\lambda_n$ includes the interval $[1, + \infty]$.

Problem 1.3 for $n=2$ and $\lambda$ in $[1/2, 1]$ was solved completely by
Beresnevich, Dickinson, Vaughan and Velani \cite{BeDiVe07,VaVe06}.

\proclaim Theorem BDVV.
For any real number $\lambda$ with $1/2 \le \lambda \le 1$, we have
$$
\dim \{ \xi \in \R : \lambda_2 (\xi) \ge \lambda\} =
\dim \{ \xi \in \R : \lambda_2 (\xi) = \lambda\} = {2 - \lambda \over 1 + \lambda}.
$$

For $n \ge 3$ the dimension of the level sets
$\{\xi \in \R : \lambda_n (\xi) = \lambda\}$ has been
determined by Schleischitz \cite{Schl16} for $\lambda > 1$.

\proclaim Theorem S.
Let $n \ge 2$ be an integer and $\lambda > 1$ a real number. Then, we have
$$
\dim \{ \xi \in \R : \lambda_n (\xi) \ge \lambda\} =
\dim \{ \xi \in \R : \lambda_n (\xi) = \lambda\} =
{2 \over n(1 + \lambda)}.    \eqno (1.4)
$$

%Theorem S extends an earlier result of Budarina, Dickinson, and Levesley \cite{BuDiLe},
%who established (1.4) for $\lambda \ge n-1$.

Let us briefly explain the easy part of the proof of Theorem S.
One way to construct a good rational approximation
$({p_1 \over q} , \ldots , {p_n \over q})$ to $(\xi, \ldots , \xi^n)$
is to start with a rational number $p/q$ very
close to $\xi$, that is, such that
$|q \xi - p| = q^{- \lambda}$, for some $\lambda > 1$.
We then observe that, for $j = 1, \ldots , n$, we have
$$
|q^j \xi^j - p^j| \ll_n q^{j-1} q^{-\lambda} \ll_n q^{n -1 -\lambda},   \quad j = 1, \ldots , n.
$$
This gives at once the lower bound
$$
\lambda_n (\xi) \ge {\lambda_1 (\xi) - n + 1 \over n},  \eqno (1.5)
$$
which is non-trivial if $\lambda_1 (\xi)$ exceeds $n$.
In particular, it then follows from (1.1) that
$$
\dim \{ \xi \in \R : \lambda_n (\xi) \ge \lambda\}
\ge \dim \{ \xi \in \R : \lambda_1 (\xi) \ge n \lambda + n - 1\}
= {2 \over n (1 + \lambda )}.     \eqno (1.6)
$$
This inequality is valid for all $\lambda$ with $\lambda \ge 1/n$, but the lower bound
is not greater than $2/ (n+1)$ for $\lambda = 1/n$, thus it is far from the truth when $n \ge 2$.
To establish Theorem S, Schleischitz proved that, for $\lambda  > 1$, all but finitely
many rational $n$-tuples which are the best
approximations of the real $n$-tuple $(\xi, \ldots , \xi^n)$
are of the form $( p/q, (p/q)^2, \ldots , (p/q)^n )$, that is, lie on the
Veronese curve $x \mapsto (x, \ldots , x^n)$.
In Section 5, we give a new, shorter (and, we believe, illuminating) proof of this
assertion.

As a first observation towards Problem 1.3 for $\lambda \le 1$,
let us note that the transference inequality
(due to Khintchine, see e.g. \cite{BuDurham})
$$
\lambda_n (\xi) \ge {w_n(\xi) \over (n-1) w_n (\xi) + n},
$$
combined with (1.2), shows that, for $1/n \le \lambda < 1/(n-1)$, we get
$$
\eqalign{
\dim \{ \xi \in \R : \lambda_n (\xi) \ge \lambda\}
& \ge  \dim \Bigl\{ \xi \in \R : w_n (\xi) \ge {n \lambda \over 1 - \lambda (n-1)} \Bigr\} \cr
& \ge {(n+1) (1 - \lambda (n-1)) \over 1 + \lambda}.  \cr }     \eqno (1.7)
$$
This (easy) lower estimate, which applies to a very small set of values of $\lambda$,
gives, unlike (1.6), that the Hausdorff dimension of the set
$\{ \xi \in \R : \lambda_n (\xi) \ge \lambda\}$ tends to $1$ as $\lambda$ tends to $1/n$.
It is superseeded by a deep result of Beresnevich \cite{Ber09}
dealing with values of $\lambda$ close to $1/n$.

\proclaim Theorem B.
Let $n \ge 2$ be an integer.
Let $\lambda$ be a real number
with $1/n \le \lambda < 3/(2n-1)$. Then, we have
$$
\dim \{ \xi \in \R : \lambda_n (\xi) \ge \lambda\} \ge
{n + 1 \over \lambda + 1} - (n-1).  \eqno (1.8)
$$

The lower bound (1.8) is not surprising since we may often expect that the
codimension of the intersection of two fractal sets is the sum of their
codimensions. Here, we intersect the Veronese curve, of dimension $1$, with
the set of real $n$-tuples $(\xi_1, \ldots , \xi_n)$ for which there exist
infinitely many integers $q$ such that
$\max_{1 \le j \le n} \| q \xi_j \| < q^{-\lambda}$,
where $\| \cdot \|$ denotes the distance to the nearest integer.
The Hausdorff dimension
of the latter set is equal to $(n+1) / (\lambda + 1)$, by a result of Dodson \cite{Do92}.

Observe that the lower bounds in (1.6) and (1.8) coincide for $\lambda = 2/n$ and are
equal to $2 / (n+2)$ at this value of $\lambda$.
Thus, it could be tempting to conjecture that we have
equalities in (1.6) and (1.8) for $\lambda \ge 2/n$ and for $\lambda \le 2/n$,
respectively. This is, however, not the case for $n \ge 3$: namely, we
show that the graph of the function $\lambda \mapsto \dim\{ \xi \in \R : \lambda_n (\xi) \ge \lambda\}$
is more complicated and presumably composed of about $n$   %%y
parts. Among our results, stated in Section 2, we extend the range of validity of (1.4)
and obtain new lower and upper bounds for the Hausdorff dimension of the set of real
numbers $\xi$ such that $\lambda_n (\xi) \ge \lambda$, for $\lambda > 1/n$.

Throughout this paper,
$\lfloor \cdot \rfloor$ denotes the integer part function
and $\lceil \cdot \rceil$ the ceiling function.
The notation $a \gg_d b$ means that $a$ exceeds $b$ times a constant
depending only on $d$. When $\gg$ is written
without any subscript, it means that the constant is absolute.
We write $a \asymp b$ if both
$a \gg b$ and $a \ll b$ hold.

\vskip 5mm

\centerline{\bf 2. Main results}

\vskip 6mm

Our first result is an extension of the range of validity of (1.4).

\proclaim Theorem 2.1.
Let $n \ge 2$ be an integer.
The spectrum of $\lambda_n$ contains the interval $[(n+4)/ (3n), + \infty]$.
Let $\lambda \ge (n+4) / (3n)$ be a real number. Then, we have
$$
\dim \{ \xi \in \R : \lambda_n (\xi) = \lambda\} =
{2 \over n(1 + \lambda)}.
$$
In particular, for any real number $\lambda$ with $\lambda > 1/3$, there exists an integer $n_0$
such that
$$
\dim \{ \xi \in \R : \lambda_n (\xi) = \lambda\} =
{2 \over n(1 + \lambda)},
$$
for any integer $n$ greater than $n_0$.

Our next result shows that the assumption `$\lambda > 1/3$' in the last
assertion of Theorem 2.1 is sharp.

\proclaim Theorem 2.2.
%For any positive real number $\lambda$ with $\lambda < 1/3$, there exists an integer $n_0$
%such that
For any integer $n \ge 2$, we have
$$
\dim \{ \xi \in \R : \lambda_n (\xi) \ge 1/3 \} \ge {2 \over (n - 1) (1 + 1/3)}.
$$
%for any integer $n$ greater than $n_0$.

Theorems 2.1 and 2.2 above are special cases of the following general statement.

\proclaim Theorem 2.3.
Let $k, n$ be integers with $1 \le k \le n$.
Let $\lambda $ be a positive real number. Then we have
$$
\dim\{ \xi \in \R : \lambda_n (\xi) \ge \lambda \} \ge
{(k+1) (1 -  (k-1) \lambda ) \over (n - k + 1) (1 + \lambda )}.      \eqno (2.1)
$$
If $\lambda > 1 / \lceil {n+1 \over 2} \rceil$, then, setting   %%y
%%DB Here should be \lceil {n+1 \over 2} \rceil becaue if n is even and h=m=(n+2)/2 we divide by zero and get infinity in (2.2)
$m = \min\{ 1 + \lfloor 1/ \lambda \rfloor , \lfloor {n+1 \over 2} \rfloor \}$, we have
$$
\dim \{ \xi \in \R : \lambda_n (\xi)  \ge \lambda\}  \le
\max_{1\le h\le m} \left\{{(h+1)(1-(h-1)\lambda) \over (n-2h+2)(1+\lambda)} \right\}.   \eqno (2.2)
$$

Observe that (1.7) is the special case $k=n$ of (2.1).
Our method does not give non-trivial upper bounds for the Hausdorff dimension
when $\lambda$ is smaller than $1 / \lfloor {n+1 \over 2} \rfloor$.

We briefly show how Theorem~2.1 follows from Theorem~2.4. For $n \ge 3$
and $\lambda \ge {n+4\over 3n}$,    %%y
we have $m \le 2$. Therefore,
for $\lambda \ge {n+4\over 3n}$, we get
$$
\max_{1\le h\le 2} \left\{{(h+1)(1-(h-1)\lambda) \over
(n-2h+2)(1+\lambda)} \right\}
= \max \left\{ {2 \over n (1 + \lambda)}, {3 (1 - \lambda) \over (n-2) (1 + \lambda)} \right\}
= {2\over n(1+\lambda)},
$$
and Theorem~2.1 follows.

As a special case of (2.1),
we obtain that, for any positive integer $m$ with $m \le n$,
$$
\dim \Bigl\{ \xi \in \R : \lambda_n (\xi) \ge {1 \over m} \Bigr\} \ge {1 \over n - m + 1}.
$$
We believe that the graph of $\lambda \mapsto \dim\{ \xi \in \R : \lambda_n (\xi) \ge \lambda\}$
is composed of about $n$ parts.

Inequality (1.5) is a special case of Lemma 3.1 of \cite{Bu10}, which asserts that,
for any positive integers $k$ and $n$ with $k$ dividing $n$, and, for any
transcendental real number $\xi$, we have
$$
\lambda_{n} (\xi) \ge {k \lambda_k (\xi) - n + k \over n}.   \eqno (2.3)
$$
Schleischitz \cite{Schl16} conjectured that (2.3) remains true when $k$ is
less than $n$ but does not divide $n$.
Our next theorem confirms this conjecture.

\proclaim Theorem 2.4. 
Let $\xi$ be a real transcendental number. 
For any positive integer $k$, we have 
$$
(k+1) \bigl( 1 + \lambda_{k+1} (\xi) \bigr) \ge k \bigl( 1 + \lambda_k (\xi) \bigr).     %%y
$$
Consequently, for every integer $n$ with $n \ge k$, we have 
$$
\lambda_n (\xi) \ge {k \lambda_k (\xi) - n + k \over n}.  
$$

Theorem 2.4 has been established independently by Schleischitz \cite{Schl19}, who also
proved a lower estimate of $\lambda_n (\xi)$ in terms of $w_k (\xi)$, for $n \ge k$.
%a stronger inequality than our first one, namely
%$$
%\lambda_n (\xi) \ge {(k-1) w_k (\xi) - k (n-k)  \over (k-1)^2 w_k (\xi) + k (n -  1)}.
%$$
%In addition, the bounds obtained in \cite{Schl19} involved uniform exponents.

The first assertion of Theorem 2.4 is of interest only when 
$\lambda_k (\xi) > 2/k$.  %$\lambda_k (\xi) > (n - k + 1)/k$. 
The last assertion is obtained by repeated 
application of the first one. This shows
at once that, if there is equality in (2.3), then we have 
$$
\lambda_m (\xi) = {k \lambda_k (\xi) - m + k \over m}, \quad m = k, \ldots , n.
$$

The present paper is organized as follows. We establish two new lower bounds for
$\lambda_n (\xi)$ in Section 3. We derive (2.1) from one of them. The second
one is Theorem 2.4 above. Section 4 is devoted to the proof of (2.2), which
follows an original approach inspired by a paper of Davenport and Schmidt \cite{DaSc69}. %%y
Finally, in Section 5, we give alternative proofs of some earlier results
of Schleischitz, including Theorem S.

\vskip 8mm

\centerline{\bf 3. Lower bounds for the exponents $\lambda_n$}

\vskip 6mm

The key ingredient for the proof of Theorem 2.5 is a new lower bound for $\lambda_n (\xi)$ in terms
of a quantity similar to $w_k(\xi)$.
For a positive integer $n$, we denote by $w_n^{\rm lead}$ the exponent of approximation defined
as in Definition 1.1, but with the
additional requirement that $|x_n|$ is not smaller than $|x_0|, \ldots , |x_{n-1}|$.
%the height of the polynomial is the absolute value of its leading coefficient

\proclaim Theorem 3.1.
Let $k, n$ be integers with $2 \le k \le n$.
Let $\xi$ be a real transcendental number. Then, we have
%$$
%\lambda_n (\xi) \ge {w_k (\xi) - k (n-k)  \over (k-1) w_k (\xi) + k (n - k + 1)},
%$$
$$
\lambda_n (\xi) \ge {w_k^{\rm lead} (\xi) - n + k \over (k-1) w_k^{\rm lead} (\xi) + n}.  \eqno (3.1)
$$
Furthermore, $w_k^{\rm lead} (\xi)$ can be replaced by $w_k (\xi)$ in (3.1) when $k=2$
or when $n = k+1$.

We suspect that $w_k^{\rm lead} (\xi)$ can be replaced by $w_k (\xi)$ in (3.1) for every $k, n$ with
$2 \le k \le n$.

\pro
Let $k, n$ be integers with $2 \le k \le n$.
Let $\xi$ be a real transcendental number.
We assume for the moment that $w_k^{\rm lead} (\xi)$ is finite and set $w_k = w_k^{\rm lead} (\xi)$.
Let $\eps$ be a positive real number.

For arbitrarily large integers $H$, there exist integers $a_0, a_1, \ldots , a_k$, not all zero,
such that $ H = |a_k| = \max\{ |a_0|, |a_1|, \ldots , |a_k|\}$ and
$$
H^{-w_k - \eps} \le \rho := |a_k \xi^k + \ldots + a_1 \xi + a_0|  \le H^{-w_k + \eps}.    \eqno (3.2)
$$
Then, by Minkowski's Theorem, there exist integers $v_0, \ldots , v_n$, not all zero, such that
$$
|v_0 \xi^j - v_j| \le |a_k|^{(n-k)/k} \rho^{1/k}, \quad 1 \le j \le k,
$$
$$
|a_0 v_i + a_1 v_{i+1} + \ldots + a_k v_{i+k}| < 1, \quad 0 \le i \le n-k.
$$
Since the $a_j$'s and $v_j$'s are integers, we get that
$$
a_0 v_i + a_1 v_{i+1} + \ldots + a_k v_{i+k} = 0, \quad 0 \le i \le n-k.
$$
Consequently,
$$
| \rho v_0 |
= | (v_0 \xi - v_1) a_1 + \ldots + (v_0 \xi^k - v_k) a_k|
\le k H |a_k|^{(n-k)/k} \rho^{1/k} = k H^{n/k} \rho^{1/k}.
$$
It then follows from (3.2) that
$$
|v_0| \le k H^{(n + (k-1) (w_k + \eps) ) / k}.      \eqno (3.3)
$$
Furthermore, for $i=1, \ldots , n-k$, we have
$$
\eqalign{
|v_0 \xi^{i+k}  - v_{i+k}|  
=  \Bigl| v_0 & 
\Bigl( {a_{k-1} \xi^{i+k-1} + \ldots + a_1 \xi^{i+1} + a_0 \xi^i \pm \rho \xi^i \over a_k} \Bigr) \cr 
& \hskip 12mm  - {a_0 v_i + a_1 v_{i+1} + \ldots + a_{k-1} v_{i+k-1} \over a_k} \Bigr|. \cr}
$$
Inductively, we derive that
$$
|v_0 \xi^{i+k} - v_{i+k}|  \ll_{n, \xi}  H^{(n-k)/k} \rho^{1/k}
 \ll_{n, \xi}  H^{(n-k-w_k + \eps) / k}, \quad i=1, \ldots , n-k.    \eqno (3.4)
$$
%Denoting by $w_k^{\rm lead}$ the exponent of approximation defined with the
%additional requirement that the height of the polynomial is the absolute
%value of its leading coefficient,
We deduce at once from (3.3) and (3.4) that
$$
\lambda_n (\xi) \ge {w_k^{\rm lead} (\xi) - n + k \over (k-1) w_k^{\rm lead} (\xi) + n}.
$$
An inspection of the proof shows that it yields $\lambda_n (\xi) \ge 1 / (k-1)$
when $w_k^{\rm lead} (\xi)$ is infinite, so (3.1) holds in all cases.

Note that the same arguments apply in the case where, instead of $|a_k|=H$,
we have $|a_k|\ge C H$, for a given positive $C$ and arbitrarily large $H$.    %%y
Also, a very similar argument shows that $w_k^{\rm lead} (\xi)$
can be replaced by $w_k^{\rm cst} (\xi)$ in (3.1), where, for $n \ge 1$, the
exponent of approximation $w_n^{\rm cst}$ is defined as in Definition 1.1,
with the additional requirement that $|x_0| \gg \max\{|x_1|, \ldots, |x_k|\}$.
Note that, for $k=2$ and $\xi\asymp 1$, equation~(3.2) implies $\max\{|a_2|, |a_0|\}\gg H$.
Consequently, for $k=2$, we can replace
$w_2^{\rm lead} (\xi)$ by $w_2 (\xi)$ in (3.1).

For the last assertion of the theorem, we need to adapt the proof of (3.1).
We proceed as follows. Let $n=k+1$, set $w_k = w_k (\xi)$ and, assuming that
$w_k$ is finite, let $\eps$ be a
positive real number.
There exists $h$ in $\{ 0, \ldots , k\}$ such that, for arbitrarily large
integers $H$, there are integers $a_0, a_1, \ldots , a_k$, not all zero, such
that $ H = \max\{ |a_0|, |a_1|, \ldots , |a_k|\}$ and
$$
H^{-w_k - \eps} \le \rho := |a_k \xi^k + \ldots + a_1 \xi + a_0|  \le H^{-w_k + \eps},
\quad |a_{h+1}\xi - a_h|\asymp H.
$$
Here and below, we set $a_{-1} = a_{k+1} = 0$.

%%DB I add this part because I think the proof here was not clear. However the proof becomes bulky. We can think if we want to add this part.
Consider the matrix 
$$
M:=\left(\matrix{
\xi&-1&0&\cdots&\cdots&\cdots&\cdots&\cdots&0\cr
\xi^2&0&-1&0&\cdots&\cdots&\cdots&\cdots&0\cr
\vdots&\vdots&\vdots&\ddots&\vdots&\vdots&&&\vdots\cr
\xi^h&0&0&\cdots&-1&0&\cdots&\cdots&0\cr
\xi^{h+2}&0&0&\cdots&0&0&-1&\cdots&0\cr
\vdots&\vdots&&&\vdots&\vdots&&\ddots&\vdots\cr
\xi^{k+1}&0&0&\cdots&0&0&\cdots&0&-1\cr
a_0&a_1&a_2&\cdots&a_h&a_{h+1}&a_{h+2}&\cdots&0\cr
0&a_0&a_1&\cdots&a_{h-1}&a_h&a_{h+1}&\cdots&a_k\cr
}\right).
$$
By expanding its determinant with respect to the
first column, we compute
$$
| \det (M) | =
\left| \sum_{j=1 \atop j\neq h+1}^{k+1} \xi^j\cdot \det\left(\matrix{a_j&a_{h+1}\cr a_{j-1}&a_h\cr} \right)
+ a_0a_h\right| = |\xi a_{h+1} - a_h| \rho =: R\rho.
$$

By Minkowski's Theorem, there exist integers $v_0, \ldots , v_n$, not all
zero, such that $v_0 > 0$,
$$
|v_0 \xi^j - v_j| \le |R|^{(n-k)/k} \rho^{1/k}, \quad 1 \le j \le n, \enspace j \not= h+1,
$$
$$
|a_0 v_0 + a_1 v_1 + \ldots + a_k v_k| < 1, \quad  |a_0 v_1 + a_1 v_2 + \ldots + a_k v_{k+1}| < 1.
$$
Since the $a_j$'s and $v_j$'s are integers, we get that
$$
a_0 v_0 + a_1 v_1 + \ldots + a_k v_k = a_0 v_1 + a_1 v_2 + \ldots + a_k v_{k+1} = 0.
$$
Observe that
$$
\sum_{i=0}^n (v_0\xi^i - v_i)(a_i\xi - a_{i-1}) = 0.
$$
Therefore,
$$
R|v_0\xi^{h+1}- v_{h+1}| \le \sum_{i=0\atop i\neq h+1}^n |a_i\xi - a_{i-1}|\cdot |v_0\xi^i - v_i|
\ll H R^{(n-k)/k}\rho^{1/k},
$$
which implies
$$
|v_0\xi^{h+1}- v_{h+1}| \ll R^{n-k\over k}\rho^{1\over k}.
$$
%$$
%|a_h| \cdot |v_0 \xi^{h+1} - v_{h+1}|
%\le | a_0 (v_1 - v_0 \xi) + a_1 (v_2 - v_0 \xi^2) + \ldots + a_{n-1} (v_n - v_0 \xi^n) |
%+ \rho |\xi| \cdot |v_0|
%$$
%and
Next, we estimate
$$
|\rho v_0|  = |a_1 (v_0 \xi - v_1) + \ldots + a_k (v_0 \xi^k - v_k)|
\ll H R^{n-k\over k} \rho^{1\over k},
$$
and the remaining part of the proof is analogous to that for~(3.1).
We omit the details.
\cqfd

\bigskip

\noi {\it Proof of the first assertion of Theorem 2.3. }
Let $k, n$ be integers with $2 \le k \le n$.
Let $\lambda > 1/n$. Inequality~(3.1)
implies that
$$
\{\xi\in \R : \lambda_n(\xi)\ge \lambda\} \subset \left\{\xi\in \R : w_k^{\rm lead}(\xi)
\ge {(\lambda+1)n-k \over 1-\lambda(k-1)}\right\}.    \eqno (3.5)
$$
Bernik \cite{Ber} established that
$$
\dim \{ \xi \in \R : w_k^{\rm lead} (\xi) \ge w \} = {k + 1 \over w + 1},    \eqno (3.6)
$$
for every real number $w$ with $w \ge k$.
The combination of (3.5) and (3.6) yields (2.1).
\cqfd

\bigskip

Similar ideas as the ones used in the proof of Theorem 3.1 allow us to bound
$\lambda_n (\xi)$ from below in terms
of $\lambda_k (\xi)$, where $k \le n$.

\medskip

\noi {\it Proof of Theorem 2.4. }
Write $\lambda_k = \lambda_k (\xi)$.
Assume that $\lambda_k$ is finite (otherwise, $\lambda_{k+1} (\xi)$ is infinite and we are done). %%y 
Let $\eps$ be a positive real number.
There exist arbitrarily large positive integers $q$ such that
$$
q^{-\lambda_k - \eps} \le \rho := \max\{ \|q \xi \|, \ldots , \| q \xi^k \| \}  \le q^{-\lambda_k + \eps}.
$$
For $j = 1, \ldots , k$, let $v_j$ be the integer such that $\|q \xi^j \| = |q \xi^j + v_j|$.
Let $h$ be in $\{1, \ldots , k\}$ such that $\rho = \| q \xi^h\|$.
We apply Minkowski's theorem to deduce the existence of integers
$a_0, a_1, \ldots , a_k$, not all zero, such that
$$
|-a_0 + a_1 \xi + \ldots + a_k \xi^k| \le \rho^{1 + ( k - 1) / (k \lambda_k)},
$$
$$
|a_j| \le \rho^{- 1 / (k \lambda_k) }, \quad 1 \le j \le k, \enspace j \not= h,
$$
$$
|a_0 q + a_1 v_1 + \ldots + a_k v_k| < 1.
$$
Since the $a_j$'s and $v_j$'s are integers, we get that
$$
a_0 q + a_1 v_1 + \ldots + a_k v_k = 0.
$$
One deduces that
$$
|-a_0 q + a_1 q \xi + \ldots + a_k q \xi^k|
= |a_1 (q \xi + v_1) + \ldots + a_k (q \xi^k + v_k)|,
$$
thus
$$
|a _h| \rho \le q \cdot \rho^{1 + ( k - 1) / (k \lambda_k)} + k \rho \cdot \rho^{- 1 / (k \lambda_k) },
$$
Since $q \le \rho^{-1 / (\lambda_k - \eps)}$, we get
$$
|a_h| \le \rho^{- c_1 \eps - 1 / (k \lambda_k) },
$$
where $c_1$, as well as $c_2, c_3, \ldots$ below, is positive and
depends at most on $n, \xi,$ and $\lambda_k$.
%{\bf [here, I have used $\lambda_k \le 1$]}
It then follows that
$$
H := \max\{ |a_0|, |a_1|, \ldots , |a_k|\}
\le  \rho^{-  c_1 \eps - 1 / (k \lambda_k) } \le  q^{c_2 \eps + 1 / k}.
$$
Using triangle inequalities as above, we get that
$$
\eqalign{
\| a_k  q \xi^{k+1} \| & \le | a_kq\xi^{k+1} - a_{k-1}v_k - a_{k-2}v_{k-1}-\ldots-a_1v_2+a_0v_1|  \cr
& \ll_{n, \xi, \eps} q   \cdot |a_k \xi^k + \ldots + a_1 \xi - a_0| +  H  \rho  \cr 
& \le q^{(1/k)  -  \lambda_k + c_3 \eps}.  \cr}   \eqno (3.7) 
$$
It follows from
$$
|a_k q| \le  q^{c_4 \eps + 1 + 1 / k} 
$$
and (3.7) that
$$
\lambda_{k+1} (\xi) \ge {\lambda_k (\xi) - 1/k \over 1 + 1/k} - c_5 \eps.
$$
As $\eps$ can be chosen arbitrarily close to $0$, we deduce that
$$
(k+1) \bigl( 1 + \lambda_{k+1} (\xi) \bigr) \ge k \bigl( 1 + \lambda_k (\xi) \bigr).
$$
This concludes the proof.
\cqfd

 \vskip 8mm

\centerline{\bf 4. Upper bound}

\vskip 6mm

Since $\lambda_n (\xi) = \lambda_n (\xi + m)$ for any integer $m$, we may assume that
$\xi\asymp 1$.
We investigate the $(n+1)$-tuples $\vp:= (q,p_1,p_2,\ldots, p_n)$ of integers
which approximate at least one point ${\bf \xi}  =
(\xi,\xi^2,\ldots, \xi^n)$ on the Veronese curve, that is, which satisfy
$$
|q\xi^i - p_i|\ll q^{-\lambda}, \quad i = 1, \ldots , n.   \eqno (4.1)
$$
Obviously, the condition $\xi\asymp 1$ is equivalent to $q\asymp p_1\asymp
p_2\asymp\cdots\asymp p_n$. For convenience, we will often write $p_0$
instead of $q$.

Throughout this section, we extensively make use of matrices of the form
$$
\Delta_{m,k}:= \left(
\matrix{
p_{k-m+1}&p_{k-m+2}&\cdots &p_k\cr
p_{k-m+2}&p_{k-m+3}&\cdots &p_{k+1}\cr
\vdots&\vdots&\ddots&\vdots\cr
p_k&p_{k+1}&\cdots&p_{k+m-1} \cr}
\right).
$$
Observe that $\Delta_{m,k}$ is an $m \times m$ matrix with $p_k$ in its antidiagonal.
Note also that the $\Delta_{m,k}$'s are precisely Hankel matrices constructed from the
sequence of $(p_k)_{k\in\{0,\ldots, n\}}$. For a given square matrix $A$ we
denote by $|A|$ the absolute value of its determinant.

%Desnanot–Jacobi identity links the determinants of different $\Delta_{m,k}$
%in the following way:
%\begin{equation}\label{eq_dj}
%|\Delta_{m+1,k}| = \frac{|\Delta_{m,k-1}||\Delta_{m,k+1}| - |\Delta_{m,k}|^2}{|\Delta_{m-1,k}|}.
%\end{equation}

\proclaim Proposition 4.1.
Assume that a tuple $\vp = (p_0,\ldots, p_n)$ in $\ZZ^{n+1}$
satisfies (4.1) for some real number $\xi$ with $\xi\asymp 1$. Then, we have
$$
|p_i\xi - p_{i+1}| \ll q^{-\lambda}, \quad \hbox{for  $i\in\{0,\ldots, n-1\}$, and}
$$
$$
|\Delta_{2,i}| \ll q^{1-\lambda}, \quad \hbox{for  $i\in\{0,\ldots, n-1\}$.}  \eqno (4.2)
$$
Conversely, if an integer tuple $\vp$ in $\ZZ^{n+1}$ with $p_0\asymp p_1\asymp\cdots
\asymp p_n$ satisfies~(4.2), then there exists a real number $\xi$ for which (4.1) is true.

\pro
For the first part of the proposition, the triangle inequality gives
$$
|p_i\xi - p_{i+1}| = |(q\xi^i+ (p_i-q\xi^i))\xi - p_{i+1}|
\le |\xi (q\xi^i - p_i)| + |q\xi^{i+1} - p_{i+1}| \ll q^{-\lambda}.
$$
For the second inequality, we have
$$
\left|\matrix{
p_{i-1}&p_i\cr
p_i&p_{i+1} \cr}
\right|
= \left|\matrix{
p_{i-1}&p_i\cr
p_i - \xi p_{i-1}&p_{i+1} - \xi p_i  \cr}
\right|
\ll q^{1-\lambda}.
$$

Finally, consider an integer tuple $\vp$ which satisfies (4.2). Then, for
$i = 1, \ldots , n-1$,  we have
$$
\left| {p_i \over  p_{i-1}} - {p_{i+1} \over p_i}\right| \ll q^{-1-\lambda}.
$$
Setting $\xi:= p_1/p_0$, these inequalities yield
$$
\left|\xi - {p_{i+1} \over p_i}\right| \ll q^{-1-\lambda}, \quad
\hbox{thus} \quad |p_i\xi - p_{i+1}|\ll q^{-\lambda}.
$$
Now we use induction on $i$. For $i=1$, the statement $|q\xi - p_1|\ll
q^{-\lambda}$ follows from the last estimate. Assuming that~(4.1)
is true for $i$, we deduce from
$$
|q\xi^{i+1} - p_{i+1}| = |(q\xi^i - p_i)\xi + p_i\xi - p_{i+1}|\ll q^{-\lambda}.
$$
that it is also true for $i+1$.
\cqfd

Proposition~4.1 allows us to investigate integer $(n+1)$-tuples~$\vp$
which satisfy~(4.2), instead of real numbers $\xi$ with
$\lambda_n(\xi)\ge \lambda$.

\proclaim Proposition 4.2. Let $\vp$ be in $\ZZ^{n+1}$ which satisfies~(4.2).
Then, for any positive integers $m,k$ with $k-m+1\ge 0$ and $k+m-1\le n$, we
have
$$
|\Delta_{m,k}|\ll q^{1 - (m-1)\lambda}.
$$

\pro By Proposition 1, there exists a real number $\xi$ which satisfies (4.1) and,
in particular, such that $|p_i\xi - p_{i+1}|\ll q^{-\lambda}$, for $ i = 1,
\ldots , n-1$. Then,
$$
|\Delta_{m,k}| = \left|\matrix{
p_{k-m+1}&\cdots &p_k\cr
p_{k-m+2}&\cdots &p_{k+1}\cr
\vdots&\ddots&\vdots\cr
p_k&\cdots&\!\!\!p_{k+m-1} \cr}
\right|
$$
is equal to
$$
\left|\matrix{
p_{k-m+1}&p_{k-m+2}&\cdots &p_k\cr
p_{k-m+2}-p_{k-m+1}\xi& p_{k-m+3}-p_{k-m+2}\xi &\cdots &p_{k+1}-p_k\xi\cr
\vdots&\vdots&\ddots&\vdots\cr
p_k-p_{k-1}\xi&p_{k+1}-p_k\xi&\cdots& p_{k+m-1}-p_{k+m-2}\xi  \cr}
\right|,
$$
which, by our assumption, is clearly $\ll q^{1-(m-1)\lambda}$.
\cqfd

The proof of Proposition 4.2 can easily be adapted to show the next proposition, which is more
general.

\proclaim Proposition 4.3.
Let $\vp$ be in $\ZZ^{n+1}$ which satisfies~(4.2) and $m$
a positive integer. For $i=0, \ldots , n-m+1$, let
$\vy_i$ denote the vector $(p_i,p_{i+1},\ldots, p_{i+m-1})$.
Then, for any sequence $c_1,c_2,\ldots, c_m$ of integers in
$\{0, \ldots , n-k+1\}$, the determinant $d(c_1,\ldots, c_m)$ of the
$m \times m$ matrix composed of the
vectors $\vy_{c_1}, \vy_{c_2},\ldots, \vy_{c_m}$ satisfies
$$
|d(c_1,\ldots, c_m)|\ll q^{1-(m-1)\lambda}.
$$

Theorem DS below is a straightforward corollary of
Theorem~3 of Davenport and Schmidt \cite{DaSc69}.

\proclaim Theorem DS.
Let $a_0,a_1,\ldots, a_h$ be integers with no common factor throughout.
Assume that, for some non-negative integers $m,k$ with $k+h-1\le m$ and $m+h\le n$, the integers
$p_k,p_{k+1},\ldots, p_{m+h}$ are related by the recurrence relation
$$
a_0p_i+a_1p_{i+1}+\cdots+ a_{h}p_{i+h} = 0, \quad k\le i\le m.
$$
Let $Z$ be the maximum of the absolute values of all the $h\times h$ determinants
formed from any $h$ of the vectors $\vy_i:=(p_i, p_{i+1}, \ldots,
p_{i+h-1})$, $i = k, \ldots , m+1$. If $Z$ is non-zero, then
$$
\max\{|a_0|,|a_1|,\ldots,|a_h|\} \ll Z^{1/(m-k-h+2)}.
$$

We are now in position to establish the second assertion of Theorem 2.3.
We use some of the ideas from~\cite{DaSc69}. Let
$\lambda> 1/\lceil(n+1)/2\rceil$  %%y
be a real number and set $m = 1 + \lfloor 1 / \lambda \rfloor$.
Let $\xi$ be a transcendental real number such that $\lambda_n (\xi) \ge \lambda$
and consider an
$(n+1)$-tuple $\vp$ for which~(4.1) is satisfied and $q$ is large enough.

Let $h$ be the smallest non-negative integer number such that the matrix
$$
P_h:=\left(\matrix{
p_{0}&p_{1}&\cdots &p_{n-h-1}&p_{n-h}\cr
p_{1}&p_{2}&\cdots &p_{n-h}&p_{n-h+1}\cr
\vdots&\vdots&\ddots&\vdots&\vdots\cr
p_h&p_{h+1}&\cdots&p_{n-1}&p_{n} \cr}
\right).
$$
has rank at most $h$. Obviously, $h\le \lceil {n+1 \over 2}\rceil$, because
for $\ell = \lceil {n+1 \over 2}\rceil$ the matrix $P_\ell$ has more rows
than columns and its rank is at most $\ell$. Also, we have $h\ge 1$ since $\vp$ is
not the zero vector. On the other hand, for $q=p_0$ large enough, we get
$h\le m$. Indeed, assume that $m<  {n+1 \over 2}$ (otherwise there is nothing
to prove) and consider $m+1$ arbitrary columns of the matrix $P_m$. By
Proposition~4.3, the matrix formed from these columns has determinant at most
$cq^{1-m\lambda}$ for some absolute positive constant $c$. Since
$\lambda>1/m$, for $q$ large enough, this determinant is zero. Since
$\lambda> 1/\lceil(n+1)/2\rceil$, we have    %%y
$$
h\le m\le \left\lfloor{n+1\over 2}\right\rfloor.  \eqno (4.3)
$$

By construction of the matrix $P_h$, there exist integers $a_0,a_1, \ldots, a_h$ with no
common factor such that
$$
a_0p_i+a_1p_{i+1}+\cdots+ a_{h}p_{i+h} = 0, \quad 0\le i\le n-h.  \eqno (4.4)
$$
Note that the matrix $P_{h-1}$ has rank $h$ and therefore the value of $Z$,
defined in Theorem~DS is non-zero. Moreover, Proposition~4.3 implies that
$Z\ll q^{1-(h-1)\lambda}$. From inequality~(4.3) we have $h-1\le n-h$ and
thus all the assumptions of Theorem~DS are satisfied. It yields
$$
H:=\max\{|a_0|,|a_1|,\ldots, |a_h|\}\le Z^{1/(n-2h+2)}\ll q^{ {1-(h-1)\lambda \over n-2h+2}}.
$$

%{\bf The next two lines are not completely clear.}
%The condition $h\le \lfloor {n+1 \over 2}\rfloor$ together with~(4.3)
%imposes an upper bound on $m$. If $n$ is even, then $m$ must be at most
%$\lfloor {n+1 \over 2}\rfloor$, otherwise $h$ may be big enough so that
%Theorem~DS is inapplicable.

Consider the relation~(4.4) for $i=0$ and divide it by $p_0=q$. Then,
the condition~(4.1) implies that
$$
|a_h\xi^h+a_{h-1}\xi^{h-1}+\ldots+a_0|\ll Hq^{-1-\lambda}\ll
H^{1 -  {(1+\lambda)(n-2h+2) \over 1 - (h-1)\lambda}}.
$$
This shows that every good
approximation $\vp$ of $\xi$ with $q$ large enough provides us with an integer polynomial
$Q_\vp (X)$ of degree at most $h$ such that $|Q_\vp(\xi)| \ll Hq^{-1-\lambda}$. Then, since
$\xi$ is transcendental, we must have infinitely many different polynomials
$Q_\vp (X)$ with this property. In other words,
$$
\{ \xi \in \RR \setminus \overline{\QQ} : \lambda_n (\xi) \ge \lambda \} \subset
\bigcup_{1\le h\le m}A_h\left(1 -  {(1+\lambda)(n-2h+2) \over 1 - (h-1)\lambda}\right),
$$
where
$$
A_h(w):= \{\xi\in\RR\;:\; |P(\xi)|\ll H(P)^{-w}\; \hbox{for infinitely many }\; P\in\ZZ[x], \deg P\le h\}.
$$
It then follows from (1.2) that
$$
\dim \{ \xi \in \RR  : \lambda_n (\xi) \ge \lambda \} \le
\max_{1\le h\le m} \left\{ {(h+1)(1-(h-1)\lambda) \over (n-2h+2)(1+\lambda)}\right\}.  
$$
The proof of the last assertion of Theorem 2.3 is complete. 

%%y Don't you think one could delete the remark?
%{\bf Remark.} For odd values of $n$, the arguments above work for arbitrarily
%small values of $\lambda$ (i.e. arbitrary large values of $m$). In that
%case,~(4.3) still imposes an upper bound $h\le (n+1)/2$ in~(4.9).
%However, if $\lambda \le {2 \over n+1}$, one can check that for $h=(n+1)/2,$ the
%expression in the maximum in~(4.9) is always bigger than 1, so we do
%not get anything better than the trivial bound $1$ for the
%Hausdorff dimension.

 \vskip 8mm

\centerline{\bf 5. A simple proof of Theorem S and further results}

\vskip 6mm

Schleischitz' proof of Theorem S (see \cite{Schl16} and Theorem 2.5.8 of \cite{BuDurham}) is
clever, but there is a simpler argument, that we present below.
The common ingredient of both proofs is the fact
that, if a rational tuple is sufficiently close to the $n$-tuple
$(\xi, \ldots , \xi^n)$, then it must lie on the Veronese curve.

Let $n \ge 2$ be an integer and $\xi$ a real number with $\lambda_n (\xi) > 1$.
Let $\lambda$ be a real number with $1 < \lambda < \lambda_n (\xi)$.
Then, there are arbitrarily large integers $q, p_1, \ldots , p_n$ such that
$$
|q \xi^j - p_j| < q^{-\lambda}, \quad j=1, \ldots , n.
$$
Set $p_0 = q$. Observe that (as in the previous section, we denote by $|A|$
the absolute value of the determinant of a square matrix $A$), for $j=1, \ldots , n-1$, we have
$$
\Delta_j := \Bigl| \matrix{p_{j-1} & p_j \cr p_j & p_{j+1} \cr}  \Bigr|
=  \Bigl| \matrix{p_{j-1}  & p_j - p_{j-1} \xi \cr p_j  & p_{j+1} - p_j \xi \cr} \Bigr|  %%y 
= |p_{j-1} (p_{j+1} - p_j  \xi) - p_j (p_j  - p_{j-1} \xi )|,   %%y 
$$
thus, by the triangle inequality,
$$
\Delta_j \ll_{\xi} |q|^{1-\lambda}.
$$
If $|q|$ is sufficiently large, then we get
$$
\Delta_1 = \ldots = \Delta_{n-1} = 0,
$$
which implies that there exist coprime non-zero integers $a, b$ such that
$$
{p_1 \over q} = {p_2 \over p_1} = \ldots = {p_n \over p_{n-1}} = {a \over b}.
$$
We deduce at once that the point
$$
\Bigl( {p_1 \over q} , \ldots , {p_n \over q} \Bigr)
= \Bigl( {a \over b} , \ldots , \Bigl( {a \over b} \Bigr)^n \Bigr)
$$
lies on the Veronese curve $x \mapsto (x, x^2, \ldots , x^n)$ and that
$q$ is an integer multiple of $b^n$. In particular, we get
$$
\Bigl| \xi - {p_1 \over q} \Bigr| = \Bigl| \xi - {a \over b} \Bigr| < q^{-1 - \lambda} \le b^{-n(1 + \lambda)}.
$$
Since $q$ (and, thus, $b$) is arbitrarily large, we deduce from the (easy half of the)
Jarn\'\i k--Besicovich theorem that
$$
\dim \{\xi \in \R : \lambda_n (\xi) \ge \lambda\} \le {2 \over n(1 + \lambda)}.
$$
Combined with (1.6), this gives a full proof of Theorem S.

%The same simple argument shows that the $p$-adic analogue of Theorem S holds.
%This improves an earlier result of Budarina et al.

\bigskip

A similar argument applies to other polynomial curves than the Veronese curve and yields
a simple proof of some results of Schleischitz \cite{Schl17}. % (and perhaps \cite{Schl18}).

\proclaim Theorem 5.1.
Let $d_1, \ldots , d_\ell$ be integers with $d_1 = 1$ and $d_1 < d_2 < \ldots < d_{\ell}$.
Let $\lambda$ be real and greater than $d_{j+1} - d_j$ for $j = 1, \ldots , \ell-1$.
Then, the Hausdorff dimension of the set of real numbers $\xi$ for which
$$
\max\{ \|q \xi^{d_1} \|, \ldots , \|q \xi^{d_\ell} \|  \} < q^{-\lambda}
$$
has infinitely many solutions in positive integers $q$, is equal to $2 / (d_\ell (\lambda + 1))$.

\pro
For $j = 2, \ldots , \ell$, let $p_j$ be an integer with $\|q \xi^{d_j} \| = |q \xi^{d_j} - p_j |$. 
Without any loss of generality, we assume that $\xi$ and $q$ are positive. %%y 
For $j = 2, \ldots , \ell-1$, we have
$$
\eqalign{
\left| \matrix{p_{j-1}^{d_{j+1} - d_j} & p_j^{d_{j+1} - d_j}
\cr p_j^{d_j - d_{j-1}}  & p_{j+1}^{d_j - d_{j-1}} \cr}  \right|
& =  \left| \matrix{p_{j-1}^{d_{j+1} - d_j}& p_j^{d_{j+1} - d_j} - p_{j-1}^{d_{j+1} - d_j} \xi^{(d_{j+1} - d_j)(d_j - d_{j-1})} \cr
p_j^{d_j - d_{j-1}} & p_{j+1}^{d_j - d_{j-1}} - p_j^{d_j - d_{j-1}}  \xi^{(d_{j+1} - d_j)(d_j - d_{j-1})} \cr} \right|  \cr
& \ll p_{j}^{d_j - d_{j-1}}  | p_{j-1}^{d_{j+1} - d_j} \xi^{(d_{j+1} - d_j)(d_j - d_{j-1})}
- p_j^{d_{j+1} - d_j} |  \cr
& \ \ \ \ + p_{j-1}^{d_{j+1} - d_j}  | p_j^{d_j - d_{j-1}}  \xi^{(d_{j+1} - d_j)(d_j - d_{j-1})}
- p_{j+1}^{d_j - d_{j-1}}  |   \cr
& \ll q^{d_j - d_{j-1}}  \cdot q^{-\lambda} + q^{d_{j+1} - d_j}  \cdot q^{-\lambda}. \cr}
$$
By our assumption, this is zero. Consequently, we get
$$
\Bigl( {p_j \over p_{j-1}} \Bigr)^{d_{j+1} - d_j} =
\Bigl( {p_{j+1} \over p_j } \Bigr)^{d_j - d_{j-1}}, \quad j = 2, \ldots , \ell-1.
$$
We deduce that
$$
\eqalign{
\Bigl( {p_j \over q} \Bigr) & = \Bigl( {p_j \over p_{j-1}} \Bigr) \times
\Bigl( {p_{j-1} \over p_{j-2}} \Bigr) \times \cdots \times \Bigl( {p_1 \over q} \Bigr) \cr
& = \Bigl( {p_{j-1} \over p_{j-2}} \Bigr)^{(d_j - d_{j-1})/(d_{j-1} - d_{j-2}) + 1}
\times \cdots \times \Bigl( {p_1 \over q} \Bigr) \cr
& = \Bigl( {p_{j-1} \over p_{j-2}} \Bigr)^{(d_j - d_{j-2})/(d_{j-1} - d_{j-2})}
\times \cdots \times \Bigl( {p_1 \over q} \Bigr) \cr
& = \Bigl( {p_{j-2} \over p_{j-3}} \Bigr)^{(d_j - d_{j-3})/(d_{j-2} - d_{j-3}) }
\times \cdots \times \Bigl( {p_1 \over q} \Bigr) \cr
& = \ldots = \Bigl( {p_2 \over p_1} \Bigr)^{(d_j - d_1)/(d_2 - d_1) } \times
\Bigl( {p_1 \over q} \Bigr) = \Bigl( {p_1 \over q} \Bigr)^{d_j}.  \cr
}
$$
This shows that there exist coprime non-zero integers $a, b$ such that
$$
\Bigl( {p_1 \over q} , \ldots , {p_n \over q} \Bigr)
= \Bigl( {a \over b} , \Bigl({a \over b}\Bigr)^{d_2},  \ldots , \Bigl( {a \over b} \Bigr)^{d_\ell} \Bigr).
$$
Since $b^{d_\ell}$ divides $q$, we conclude as above.

%the best rational approximants are lying on the
%curve $\xi \mapsto (\xi, \xi^{d_2}, \ldots , \xi^{d_\ell})$.

\bigskip

Similar arguments allow us to give an alternative proof of a result of Schleischitz
asserting that the inequality
$$
\hla_n (\xi) \le \max \Bigl\{ {1 \over n}, {1 \over \lambda_1 (\xi)} \Bigr\}     \eqno (5.1)
$$
holds, where $\hla_n (\xi)$ is the supremum of the real numbers $\lambda$ 
for which the inequalities (1.3)
have a non-zero integer solution for all sufficiently large $X$.
We do not claim that the proof below is simpler than the original one.

Since (5.1) is clearly true for $n=1$ and for $\lambda_1 (\xi) = 1$,
we assume that $n \ge 2$ and $\lambda_1 (\xi) > 1$.
Let $q$ be a large positive integer and $v$ be a real number greater than $1$ such that
$$
\| q \xi \| = | q \xi - p| \le q^{-v}.
$$
Then, we check that
$$
| q^j \xi^j - p^j | \ll q^{j-1-v}, \quad 1 \le j \le n.
$$
Let $v'$ be a real number with $1 < v' < \min\{v, n\}$ and set $X = q^{v'}$.
Let $x$ be a positive integer with $x < X$.
Assume that $\hla_n (\xi) > 1 / v'$.
Then, there are integers $x_1, \ldots , x_n$ such that
$$
|x \xi^j - x_j| \ll X^{-1 / v'}.
$$
We have
$$
\Bigl| \matrix{q & p \cr x & x_1 \cr}  \Bigr|
= \Bigl| \matrix{q & q \xi -p \cr x& x \xi - x_1\cr}  \Bigr| \ll X q^{-v} + q X^{-1 / v'} < 1,
$$
if $q$ is large enough.
As $\gcd (p, q) = 1$, we derive that $q$ divides $x$. Thus, the determinant
$$
\Bigl| \matrix{q^2 & p^2 \cr x & x_2 \cr}  \Bigr|
$$
is an integer multiple of $q$. However, it satisfies
$$
\Bigl| \matrix{q^2 & p^2 \cr x & x_2 \cr}  \Bigr|
= \Bigl| \matrix{q^2 & q^2 \xi^2 -p^2 \cr x& x \xi^2 - x_2 \cr}  \Bigr| \ll X q^{1-v} + q^2 X^{-1 / v'} < q.
$$
Consequently, we derive that, if $q$ is large enough,
the determinant is equal to $0$, hence, $q^2$ divides $x$.
Continuing in the same way, we deduce that $q^n$ divides $x$, a contradiction
with the inequalities $1 \le x < q^n$.

\bigskip
\noindent
{\bf Acknowledgement}.
The main part of this work has been done 
while Yann Bugeaud was visiting the University of Sydney,
supported by the Sydney Mathematical Research Institute International Visitor Program.

\vskip 8mm

\goodbreak

\centerline{\bf References}

\vskip 5mm

\beginthebibliography{999}

\bibitem{Ber09}
V. Beresnevich,
{\it Rational points near manifolds and metric
Diophantine approximation},
 Ann. of Math. (2) 175 (2012),  187--235.

\bibitem{BeDiVe07}
V. Beresnevich, D. Dickinson and S. L. Velani,
{\it Diophantine approximation on planer curves
and the distribution of rational points}, with
an appendix by R.C. Vaughan: `Sums of two squares near perfect squares',
Ann. of Math. 166 (2007), 367--426.

\bibitem{Ber}
V. I. Bernik,
{\it Application of the Hausdorff dimension in the theory of
Diophantine approximations}, Acta Arith. 42 (1983), 219--253 (in
Russian).
English transl. in Amer. Math. Soc. Transl. 140 (1988), 15--44.

%\bibitem{BuDiLe}
%N. Budarina, D. Dickinson, and J. Levesley,
%{\it Simultaneous Diophantine approximation on
%polynomial curves}.
%Mathematika 56 (2010), 77--85.
\bibitem{BuLiv}
Y. Bugeaud,
Approximation by algebraic numbers,
Cambridge Tracts in Mathematics, Cambridge, 2004.

\bibitem{Bu10}
Y. Bugeaud,
{\it On simultaneous rational approximation
to a real number and its integral powers},
Ann. Inst. Fourier (Grenoble) 60 (2010), 2165--2182.

\bibitem{BuFARF}
Y. Bugeaud,
{\it Hausdorff dimension and Diophantine approximation.}
In: Further Developments in Fractals and Related Fields, J. Barral, S. Seuret (Eds),
Birkh\"auser, 2013, pp. 35-45.

\bibitem{BuDurham}
Y. Bugeaud,
{\it Exponents of Diophantine approximation}.
In: Dynamics and Analytic Number Theory.
Proceedings of the Durham Easter School 2014.
Edited by D. Badziahin, A. Gorodnik, N. Peyerimhoff.

\bibitem{BuLa05a}
Y. Bugeaud and M. Laurent,
{\it Exponents of Diophantine Approximation and
Sturmian Continued Fractions},
Ann. Inst. Fourier (Grenoble) 55 (2005), 773--804.

\bibitem{DaSc69}
{H. Davenport and W. M. Schmidt},
{\it Approximation to real numbers by
algebraic integers}, Acta Arith. {15} (1969), 393--416.

\bibitem{Do92}
M. M. Dodson,
{\it Hausdorff dimension, lower order and Khintchine's
theorem in metric Diophantine approximation},
J. reine angew. Math. 432 (1992), 69--76.

\bibitem{Jar31}
V. Jarn\'\i k,
{\it \"Uber die simultanen Diophantische
Approximationen}, Math. Z. 33 (1931), 505--543.

\bibitem{Mah32}
K. Mahler,
{\it Zur Approximation der Exponentialfunktionen und des
Logarithmus. I, II},
J. reine angew. Math. 166 (1932), 118--150.

\bibitem{Schl16}
J. Schleischitz,
{\it On the spectrum of Diophantine approximation constants},
Mathematika 62 (2016), 79--100.

\bibitem{Schl17}
J. Schleischitz,
{\it Diophantine approximation on polynomial curves },
Math. Proc. Cambridge Philos. Soc. 163 (2017), 533--546.

\bibitem{Schl19}
J. Schleischitz,
{\it Going-up theorems for simultaneous Diophantine approximation}.
Preprint.

\bibitem{VaVe06}
R. C. Vaughan and S. Velani,
{\it Diophantine approximation on planar curves: the convergence theory},
Invent. Math. 166 (2006), 103--124.

\endthebibliography

\vskip1cm

\noindent Yann Bugeaud

\noindent Universit\'e de Strasbourg

\noindent Math\'ematiques

\noindent 7, rue Ren\'e Descartes

\noindent 67084 STRASBOURG  (FRANCE)

\vskip2mm

\noindent {\tt bugeaud@math.unistra.fr}

\vskip1cm

\noindent Dzmitry Badziahin

\noindent The University of Sydney

\noindent Camperdown 2006, NSW (Australia)

\vskip2mm

\noindent {\tt dzmitry.badziahin@sydney.edu.au}

\bye